# EXISTENCE, UNIQUENESS AND APPROXIMATION OF A STOCHASTIC SCHRÖDINGER EQUATION: THE DIFFUSIVE CASE

By Clément Pellegrini

*Université Lyon 1, Institut Camille Jordan*

Recent developments in quantum physics make heavy use of so-called "quantum trajectories." Mathematically, this theory gives rise to "stochastic Schrödinger equations," that is, perturbation of Schrödinger-type equations under the form of stochastic differential equations. But such equations are in general not of the usual type as considered in the literature. They pose a serious problem in terms of justifying the existence and uniqueness of a solution, justifying the physical pertinence of the equations. In this article we concentrate on a particular case: the diffusive case, for a two-level system. We prove existence and uniqueness of the associated stochastic Schrödinger equation. We physically justify the equations by proving that they are a continuous-time limit of a concrete physical procedure for obtaining a quantum trajectory.

**1. Introduction.** Belavkin equations (also called stochastic Schrödinger equations) are classical stochastic differential equations describing the evolution of an open quantum system undergoing a continuous quantum measurement. The solutions of such equations are called quantum trajectories and describe the time evolution of the state of the system. The random nature of the result of quantum measurement is at the origin of the stochastic character of the evolution.

The first rigorous description of a state undergoing a continuous measurement is due to Davies in [4]. It describes in quantum optics, the behavior of an atom from which we observe the photon emission. This is the so-called "resonance fluorescence" experiment (see [6] and [2]).

In the literature, essentially two kinds of Belavkin equations are considered: they are driven either by a Brownian motion or by a counting process.

But the kind of equations which are obtained this way are of nonusual type compared to the usual theory of stochastic differential equations. In particular there is no reference in physics nor in mathematics, where the existence and the uniqueness of the solution of such equations are discussed. Furthermore, the physical justification of the apparition of these equations requires in general a quite heavy mathematical framework (Von Neumann algebra, conditional expectation, filtering, etc.). The high technology of such tools contrasts with the simplicity and the intuition of the physical model.

An approach to such equations, which is physically very intuitive, is the one of repeated quantum interactions. The setup is the following. The continuous measurement model is obtained as a limit of discrete models. This discrete model is a naive approach to the interaction of a simple system interacting with a field. The field is represented as a chain of independent copies of small pieces of environment. The simple system interacts, for a time interval $h$, with one piece of the environment. After that interaction an observable of the piece of environment is measured. The random result of the measurement induces a random new state for the small system. The small system then interacts again with another piece of the environment for a time interval $h$. A measurement of the same observable of this second copy is performed. And so on.

This experiment gives rise to a discrete evolution of the state of the small system, which is a Markov chain. The continuous-time limit ($h \to 0$) of this evolution should give rise to the quantum trajectories.

Repeated quantum interactions have been considered by Attal and Pautrat in [1] and by Gough in [5]. The continuous limit of repeated quantum interactions is rigorously shown to converge to a quantum stochastic differential equation in [1]. The setup of measuring an observable of the chain after each interaction is considered in [5], but the continuous limit, the existence and the uniqueness of the solutions are not all treated rigorously in this reference.

The aim of this article is to study the diffusive Belavkin equation, to show existence and uniqueness of the solution, to show its approximation by repeated quantum interaction models. The same results for the equation concerning the counting process are developed in another article [11].

This article is structured as follows: In Section 1, we present the mathematical model of repeated quantum interactions with measurement. We define discrete-time quantum trajectories and focus on their probabilistic properties. In particular, it is shown that these processes are classical Markov chains. Finally we deal with the model of a two-level atom in interaction with a spin chain and we describe the discrete stochastic evolution equations in this setting.

Section 2 is devoted to the continuous model. We present the two different types of Belavkin equations whose solutions are continuous-time quantum trajectories. We then prove existence and uniqueness of solutions in the diffusive case.

The link between discrete and continuous models is provided in Section 3. It is shown that particular discrete quantum trajectories (for a two-level model) satisfy stochastic equations which are discrete-time diffusive equations. We use result of weak convergence of stochastic integrals in order to prove that solutions of diffusive Belavkin stochastic equations are obtained as a limit of discrete trajectories.

**2. Discrete quantum trajectories.** We make here precise the mathematical framework to describe the model of discrete quantum trajectories.

2.1. *Repeated quantum measurements.* The physical setup is the one of a small quantum system, represented by a Hilbert space $\mathcal{H}_0$, coupled to a field modeled by an infinite chain of identical independent quantum systems. Each piece of the field is represented by a Hilbert space $\mathcal{H}$. Each space is equipped with a positive, trace-class operator with trace 1. This operator is called a "state" or "a density matrix." In this section, we present the random character of repeated measurements.

The discrete model of interaction is called "quantum repeated interactions." Each copy $\mathcal{H}$ of the environment interacts with $\mathcal{H}_0$, one after the other, during a time interval of length $h$. Information on the evolution of the small system is obtained by performing a measurement of $\mathcal{H}$ after each interaction.

For the first interaction, the compound system is described by the tensor product $\mathcal{H}_0 \otimes \mathcal{H}$ and the interaction is characterized by a total Hamiltonian $H_{\text{tot}}$ which is a self-adjoint operator on $\mathcal{H}_0 \otimes \mathcal{H}$. Its general form is

$$(2.1) \qquad H_{\text{tot}} = H_0 \otimes I + I \otimes H + H_{\text{int}},$$

where $H_0$ and $H$ are the free Hamiltonians of each system and $H_{\text{int}}$ is the interaction Hamiltonian. The operator

$$U = e^{ihH_{\text{tot}}}$$

describes the first interaction as follows. In the Schrödinger picture, if $\rho$ denotes any state on the tensor product, the evolution is given by

$$\rho \mapsto U\rho U^\star.$$

After this first interaction, a second copy of $\mathcal{H}$ interacts with $\mathcal{H}_0$ in the same way. And so on.

As the field is supposed to be an infinite chain, the whole sequence of successive interactions is described by the state space

$$(2.2) \qquad \mathbf{\Gamma} = \mathcal{H}_0 \otimes \bigotimes_{k \geq 1} \mathcal{H}_k,$$

where $\mathcal{H}_k$ designs the $k$th copy of $\mathcal{H}$. The countable tensor product $\bigotimes_{k\geq 1}\mathcal{H}_k$ is defined as follows. We consider that $\mathcal{H}$ is finite dimensional and that $\{X_0, X_1, \ldots, X_n\}$ is a fixed orthonormal basis of $\mathcal{H}$. The orthogonal projector onto $\mathbb{C}X_0$ is denoted by $|X_0\rangle\langle X_0|$ (this is the braket notation in mathematical physics; see the remark below). This is the ground state (or vacuum state) of $\mathcal{H}$. The tensor product is taken with respect to $X_0$ (for details, see [1]), that is, we define an orthonormal basis of $\bigotimes_{k\geq 1}\mathcal{H}_k$ with respect to this vector. It is described as follows.

Let $\mathcal{P}$ be the set of finite subset $A$ of the form $A = \{(n_1, i_1), \ldots, (n_k, i_k)\}$ of $\mathbb{N}^\star \times \{1, \ldots, n\}$ such that the $n_i$'s are two-by-two disjoint. The orthonormal basis of $\bigotimes_{k\geq 1}\mathcal{H}_k$ with respect to $X_0$ is the family

$$\{X_A, A \in \mathcal{P}\},$$

where for $A = \{(n_1, i_1), \ldots, (n_k, i_k)\}$, we define $X_A$ as the vector

$$X_0 \otimes \cdots \otimes X_{i_1} \otimes X_0 \otimes \cdots \otimes X_0 \otimes X_{i_2} \otimes \cdots,$$

of $\bigotimes_{k\geq 1}\mathcal{H}_k$, where $X_{i_j}$ appears in the copy number $n_j$ of $\mathcal{H}$. The infinite tensor product allows us to work in a single space but the structure of Hilbert space does not appear explicitly in the rest of the paper.

REMARK. A vector $Y$ in a Hilbert space $\mathcal{H}$ is represented by the application $|Y\rangle$ from $\mathbb{C}$ to $\mathcal{H}$ which acts with the following way: $|Y\rangle(\lambda) = |\lambda Y\rangle$. The linear form on $\mathcal{H}$ is represented by the operators $\langle Z|$ which act on the vector $|Y\rangle$ by $\langle Z||Y\rangle = \langle Z, Y\rangle$, where $\langle \cdot, \cdot \rangle$ denotes the scalar product of $\mathcal{H}$.

The unitary evolution describing the $k$th interaction is given by the operator $U_k$ which acts as $U$ on $\mathcal{H}_0 \otimes \mathcal{H}_k$, whereas it acts as the identity operator on the other copies of $\mathcal{H}$. If $\rho$ is a state on $\mathbf{\Gamma}$, the effect of the $k$th interaction is

$$\rho \mapsto U_k \rho U_k^\star.$$

Hence the result of the $k$ first interactions is described by the operator $V_k$ on $\mathcal{B}(\mathbf{\Gamma})$ defined by the recursive formula

$$(2.3) \qquad \begin{cases} V_{k+1} = U_{k+1}V_k, \\ V_0 = I, \end{cases}$$

and the evolution of states is then given, in the Schrödinger picture, by

$$(2.4) \qquad \rho \mapsto V_k \rho V_k^\star.$$

We present now the indirect measurement principle. The idea is to perform a measurement of an observable of the field after each interaction.

A measurement of an observable of $\mathcal{H}_k$ is modeled as follows. Let $A$ be any observable on $\mathcal{H}$, with spectral decomposition $A = \sum_{j=1}^{p} \lambda_j P_j$. We consider its natural ampliation on $\boldsymbol{\Gamma}$:

$$(2.5) \qquad A^k := \bigotimes_{j=0}^{k-1} I \otimes A \otimes \bigotimes_{j \geq k+1} I.$$

The result of the measurement of $A^k$ is random; the accessible data are its eigenvalues. If $\rho$ denotes the reference state of $\boldsymbol{\Gamma}$, the observation of $\lambda_i$ is obtained with probability

$$P[\text{to observe } \lambda_j] = \text{Tr}[\rho P_j^k], \qquad j \in \{1, \ldots, p\},$$

where $P_i^k$ is the ampliation of $P_i$ in the same way as (2.5). If we have observed the eigenvalue $\lambda_j$, the "projection postulate" called "wave packet reduction" imposes that the state after measurement is

$$\rho_j = \frac{P_j^k \rho P_j^k}{\text{Tr}[\rho P_j^k]}.$$

REMARK. This corresponds to the new reference state depending on the result of the observation. Another measurement of the observable $A^k$ (with respect to this new state) would give $P[\text{to observe } \lambda_j] = 1$ (for $P_i P_j = 0$ if $i \neq j$). This means that only one measurement after each interaction gives significant information. This justifies the principle of repeated interactions.

The repeated quantum measurements are the combination of the previous description and the successive interactions (2.4). After each interaction, the measurement procedure involves a random modification of the system. It defines namely a sequence of random states which is called "discrete quantum trajectory."

The initial state on $\boldsymbol{\Gamma}$ is chosen to be

$$\mu = \rho \otimes \bigotimes_{j \geq 1} \beta_j$$

where $\rho$ is any state on $\mathcal{H}_0$ and each $\beta_i = \beta$ is a fixed state on $\mathcal{H}$. We denote by $\mu_k$ the new state after $k$ interactions, that is,

$$\mu_k = V_k \mu V_k^\star.$$

The probability space describing the experience of repeated measurements is $\Omega^{\mathbb{N}^\star}$ where $\Omega = \{1, \ldots, p\}$. The integers $i$ correspond to the indexes of the eigenvalues of $A$. We endow $\Omega^{\mathbb{N}^\star}$ with the cylinder $\sigma$-algebra generated by the sets:

$$\Lambda_{i_1, \ldots, i_k} = \{\omega \in \Omega^{\mathbb{N}^\star} / \omega_1 = i_1, \ldots, \omega_k = i_k\}.$$

Note that for all $j$, the unitary operator $U_j$ commutes with all $P^k$, for $k < j$. For any set $\{i_1, \ldots, i_k\}$, we can define the following nonnormalized state:

$$\tilde{\mu}(i_1, \ldots, i_k) = (I \otimes P_{i_1} \otimes \cdots \otimes P_{i_k} \otimes I \cdots) \mu_k (I \otimes P_{i_1} \otimes \cdots \otimes P_{i_k} \otimes I \cdots)$$
$$= (P_{i_k}^k \cdots P_{i_1}^1) \mu_k (P_{i_1}^1 \cdots P_{i_k}^k).$$

It is the nonnormalized state which corresponds to the successive observation of the eigenvalues $\lambda_{i_1}, \ldots, \lambda_{i_k}$ during the $k$ first measurements. The probability to observe these eigenvalues is

$$P[\Lambda_{i_1, \ldots, i_k}] = P[\text{to observe } (\lambda_{i_1}, \ldots, \lambda_{i_k})] = \text{Tr}[\tilde{\mu}(i_1, \ldots, i_k)].$$

This way, we define a probability measure on the cylinder sets of $\Omega^{\mathbb{N}^\star}$ which satisfies the Kolmogorov Consistency Criterion. Hence it defines a unique probability measure on $\Omega^{\mathbb{N}^\star}$. The discrete quantum trajectory on $\boldsymbol{\Gamma}$ is then given by the following random sequence of states:

$$\tilde{\rho}_k : \Omega^{\mathbb{N}^\star} \longrightarrow \mathcal{B}(\boldsymbol{\Gamma}),$$
$$\omega \longmapsto \tilde{\rho}_k(\omega_1, \ldots, \omega_k) = \frac{\tilde{\mu}(\omega_1, \ldots, \omega_k)}{\text{Tr}[\tilde{\mu}(\omega_1, \ldots, \omega_k)]}.$$

From the construction and the remarks above, the following is immediate.

PROPOSITION 2.1. *Let $(\tilde{\rho}_k)$ be the above random sequence of states; we have for all $\omega \in \Omega^{\mathbb{N}^\star}$:*

$$\tilde{\rho}_{k+1}(\omega) = \frac{P_{\omega_{k+1}}^{k+1} U_{k+1} \tilde{\rho}_k(\omega) U_{k+1}^\star P_{\omega_{k+1}}^{k+1}}{\text{Tr}[\tilde{\rho}_k(\omega) U_{k+1}^\star P_{\omega_{k+1}}^{k+1} U_{k+1}]}.$$

The following theorem is an easy consequence of the previous proposition.

THEOREM 2.1. *The discrete quantum trajectory $(\tilde{\rho}_n)_n$ is a Markov chain, with values on the set of states of $\mathcal{H}_0 \bigotimes_{i \geq 1} \mathcal{H}_i$. It is described as follows:*

$$P[\tilde{\rho}_{n+1} = \mu / \tilde{\rho}_n = \theta_n, \ldots, \tilde{\rho}_0 = \theta_0] = P[\tilde{\rho}_{n+1} = \mu / \tilde{\rho}_n = \theta_n].$$

*If $\tilde{\rho}_n = \theta_n$, the random state $\tilde{\rho}_{n+1}$ takes one of the values:*

$$\frac{P_i^{n+1}(U_{n+1}(\theta_n \otimes \beta) U_{n+1}^\star) P_i^{n+1}}{\text{Tr}[(U_{n+1}(\theta_n \otimes \beta) U_{n+1}^\star) P_i^{n+1}]}, \qquad i = 1, \ldots, p,$$

*with probability* $\text{Tr}[(U_{n+1}(\theta_n \otimes \beta) U_{n+1}^\star) P_i^{n+1}]$.

In general, one is more interested in the reduced state on the small system $\mathcal{H}_0$ only. This state is given by taking a partial trace on $\mathcal{H}_0$. Let us recall what partial trace is. If $\mathcal{H}$ is any Hilbert space, we denote by $\text{Tr}_\mathcal{H}[W]$ the trace of a trace-class operator $W$ on $\mathcal{H}$.

DEFINITION 2.1. Let $\mathcal{H}$ and $\mathcal{K}$ be two Hilbert spaces. If $\alpha$ is a state on a tensor product $\mathcal{H} \otimes \mathcal{K}$, then there exists a unique state $\eta$ on $\mathcal{H}$ which is characterized by the property

$$\mathrm{Tr}_{\mathcal{H}}[\eta X] = \mathrm{Tr}_{\mathcal{H} \otimes \mathcal{K}}[\alpha(X \otimes I)]$$

for all $X \in \mathcal{B}(\mathcal{H})$. This unique state $\eta$ is called the partial trace of $\alpha$ on $\mathcal{H}$ with respect to $\mathcal{K}$.

Let $\mathbf{E}_0(\alpha)$ denote the partial trace on $\mathcal{H}_0$ with respect to $\bigotimes_{k \geq 1} \mathcal{H}_k$ of any state $\alpha$ on $\boldsymbol{\Gamma}$. We define a random sequence of states on $\mathcal{H}_0$ as follows. For all $\omega$ in $\Omega^{\mathbb{N}^\star}$, define the discrete quantum trajectory on $\mathcal{H}_0$

(2.6) $$\rho_n(\omega) = \mathbf{E}_0[\tilde{\rho}_n(\omega)].$$

An immediate consequence of Theorem 2.1 is the following result.

THEOREM 2.2. *The quantum trajectory $(\rho_n)_n$ defined by formula (2.6) is a Markov chain with values in the set of states on $\mathcal{H}_0$. If $\rho_n = \chi_n$, then $\rho_{n+1}$ takes one of the values:*

$$\mathbf{E}_0\left[\frac{(I \otimes P_i) U(\chi_n \otimes \beta) U^\star (I \otimes P_i)}{\mathrm{Tr}[U(\chi_n \otimes \beta) U^\star (I \otimes P_i)]}\right], \qquad i = 1, \ldots, p,$$

*with probability* $\mathrm{Tr}[U(\chi_n \otimes \beta) U^\star (I \otimes P_i)]$.

REMARK 1. Let us stress that

$$\frac{(I \otimes P_i) U(\chi_n \otimes \beta) U^\star (I \otimes P_i)}{\mathrm{Tr}[U(\chi_n \otimes \beta) U^\star (I \otimes P_i)]}$$

is a state on $\mathcal{H}_0 \otimes \mathcal{H}$. In this situation, the notation $\mathbf{E}_0$ denotes the partial trace on $\mathcal{H}_0$ with respect to $\mathcal{H}$. The infinite tensor product $\boldsymbol{\Gamma}$ is just needed to have a clear description of the repeated interactions and the probability space $\Omega^{\mathbb{N}^\star}$.

The next section is devoted to the particular case of a two-level atom in contact with a photon stream. Because of physical considerations, this case is often the central case in the literature concerning continuous measurement.

2.2. *A two-level atom.* The Hilbert spaces describing the physical situation are now $\mathcal{H}_0 = \mathcal{H} = \mathbb{C}^2$.

In this section, we establish a discrete quantum evolution equation for $(\rho_n)$ which is a discrete approximation of the Belavkin equation.

The main goal of this section is to obtain a formula of the following form:

(2.7) $$\rho_{k+1} = f(\rho_k, X_{k+1}),$$

where $(X_k)_k$ is a sequence of random variables. Such a formula is obtained from the description of Theorem 2.2 and the computation of the partial trace operation.

The state $\rho_k$ can be considered as an initial state (according to the Markov property of Theorem 2.2). Thus we can consider a single interaction with a system $(\mathcal{H}, \beta)$ [actually this is the $(k+1)$st copy]. We consider an observable of the form $A = \lambda_0 P_0 + \lambda_1 P_1$ and the unitary operator describing the interaction is a unitary $4 \times 4$ unitary matrix.

In order to compute the state given by the projection postulate and the partial trace, we choose a suitable basis. If $(X_0 = \Omega, X_1 = X)$ is an orthonormal basis of $\mathbb{C}^2$, for the space $\mathcal{H}_0 \otimes \mathcal{H}$ we consider the following basis: $\Omega \otimes \Omega, X \otimes \Omega, \Omega \otimes X, X \otimes X$. In this basis, the unitary operator $U$ can be written by blocks in the following way:

$$U = \begin{pmatrix} L_{00} & L_{01} \\ L_{10} & L_{11} \end{pmatrix}$$

where each $L_{ij}$ is an operator on $\mathcal{H}_0$. For $\beta$ we choose

$$\beta = |\Omega\rangle\langle\Omega|.$$

As a consequence, the state after the interaction is

$$(2.8) \qquad \mu_{k+1} = U(\rho_k \otimes \beta)U^\star = \begin{pmatrix} L_{00}\rho_k L_{00}^\star & L_{00}\rho_k L_{10}^\star \\ L_{10}\rho_k L_{00}^\star & L_{10}\rho_k L_{10}^\star \end{pmatrix}.$$

Thanks to the description of Theorem 2.2, we define the two possible nonnormalized states

$$(2.9) \qquad \mathcal{L}_0(\rho_k) = \mathbf{E}_0[I \otimes P_0(\mu_{k+1})I \otimes P_0],$$

$$(2.10) \qquad \mathcal{L}_1(\rho_k) = \mathbf{E}_0[I \otimes P_1(\mu_{k+1})I \otimes P_1].$$

These are operators on $\mathcal{H}_0$; the nonnormalized state $\mathcal{L}_0(\rho_k)$ appears with probability $p_{k+1} = \mathrm{Tr}[\mathcal{L}_0(\rho_k)]$ and $\mathcal{L}_1(\rho_k)$ with probability $q_{k+1} = \mathrm{Tr}[\mathcal{L}_1(\rho_k)]$.

Let us define the random variable $\nu_{k+1}$ on $\{0, 1\}$:

$$\begin{cases} \nu_{k+1}(0) = 0, & \text{with probability } p_{k+1}, \\ \nu_{k+1}(1) = 1, & \text{with probability } q_{k+1}. \end{cases}$$

With these notation, the discrete quantum trajectory can be described as follows. For all $\omega \in \Omega^{\mathbb{N}^\star}$, we have

$$(2.11) \qquad \rho_{k+1}(\omega) = \frac{\mathcal{L}_0(\rho_k(\omega))}{p_{k+1}(\omega)}(1 - \nu_{k+1}(\omega)) + \frac{\mathcal{L}_1(\rho_k(\omega))}{q_{k+1}(\omega)}\nu_{k+1}(\omega).$$

In order to obtain the final discrete quantum evolution equation, we consider the centered and normalized random variable

$$X_{k+1} = \frac{\nu_{k+1} - q_{k+1}}{\sqrt{q_{k+1}p_{k+1}}}.$$

We define the associated filtration on $\{0,1\}^{\mathbf{N}^\star}$:
$$\mathcal{F}_k = \sigma(X_i, i \leq k).$$
By construction, we have $\mathbf{E}[X_{k+1}/\mathcal{F}_k] = 0$ and $\mathbf{E}[X_{k+1}^2/\mathcal{F}_k] = 1$. Thus we can write the discrete evolution equation for the quantum trajectory in terms of the random variables $(X_k)$:

$$(2.12) \quad \rho_{k+1} = \mathcal{L}_0(\rho_k) + \mathcal{L}_1(\rho_k) + \left[ -\sqrt{\frac{q_{k+1}}{p_{k+1}}} \mathcal{L}_0(\rho_k) + \sqrt{\frac{p_{k+1}}{q_{k+1}}} \mathcal{L}_1(\rho_k) \right] X_{k+1}.$$

The above equation can be considered in a general way and the unique solution starting from $\rho_0$ is the discrete quantum trajectory described in Theorem 2.2. Let us stress that this sequence depends on the length of time of interaction. This dependence will allow us to prove a continuous-time approximation result in Section 3. For the moment, the next section is devoted to describing continuous-time quantum trajectories.

**3. Belavkin equations.** As was announced in the Introduction, it is commonly assumed that the evolution of a system undergoing a continuous measurement is described by stochastic differential equations. A model of interaction can be provided to describe an atom in contact with a continuous field. In this setting, the description of the principle of indirect measurement needs highly technical tools in order to obtain rigorous statements. Such theories are not the purpose of this article. We just give the physical setup in order to introduce the Belavkin stochastic differential equations.

Consider a two-level system, described by $\mathbb{C}^2$, in interaction with an environment (classically described by a Fock space). The time evolution is given by a unitary process $(U_t)$ which satisfies a quantum stochastic differential equation (cf. [10]). Without measurement the evolution of the small system is given by a norm-continuous semigroup $\{T_t\}_{t \geq 0}$. The Linblad generator of $(T_t)$ is denoted by $L$ and we have the "Master equation":
$$\frac{d\rho_t}{dt} = L(\rho_t) = -i[H_0, \rho_t] - \frac{1}{2}\{CC^\star, \rho_t\} + C\rho_t C^\star,$$
where $C$ is any operator on $\mathbb{C}^2$ and $H_0$ is the Hamiltonian of the atom.

In the theory of time-continuous measurement $L$ is decomposed as the sum of $\mathcal{L} + \mathcal{J}$ where $\mathcal{J}$ represents the instantaneous state change taking place when detecting a photon, and $\mathcal{L}$ describes the smooth state variation in between these instants. These operators are defined by
$$\mathcal{L}(\rho) = -i[H_0, \rho] - \tfrac{1}{2}\{CC^\star, \rho\},$$
$$\mathcal{J}(\rho) = C\rho C^\star.$$

Thanks to the work of Davies in [4], two types of stochastic differential equations can be derived. The solutions of these equations are then called "continuous-time quantum trajectories":

1. The "diffusive equation" (homodyne detection experiment) is given by

$$(3.1) \qquad d\rho_t = L(\rho_t)\,dt + [\rho_t C^\star + C\rho_t - \text{Tr}(\rho_t(C + C^\star))\rho_t]\,dW_t,$$

where $W_t$ describes a one-dimensional Brownian motion.

2. The "jump equation" (resonance fluorescence experiment) is

$$(3.2) \qquad d\rho_t = L(\rho_t)\,dt + \left[\frac{\mathcal{J}(\rho_t)}{\text{Tr}[\mathcal{J}(\rho_t)]} - \rho_t\right](d\tilde{N}_t - \text{Tr}[\mathcal{J}(\rho_t)]\,dt),$$

where $\tilde{N}_t$ is a counting process with stochastic intensity $\int_0^t \text{Tr}[\mathcal{J}(\rho_s)]\,ds$.

A physical justification of (3.1) as limit of discrete quantum trajectories is given in Section 3. For the moment, we shall focus on the general problem of existence and uniqueness of a solution of (3.1). The jump equation and all the convergence theorems referring to this case are treated in detail in [11] with different techniques.

3.1. *Existence and uniqueness.* Let $\rho_0$ be any state; we aim to show existence and uniqueness for the stochastic differential equation

$$(3.3) \quad \rho_t = \rho_0 + \int_0^t L(\rho_s)\,ds + \int_0^t [\rho_s C^\star + C\rho_s - \text{Tr}(\rho_s(C + C^\star))\rho_s]\,dW_s.$$

Classical theorems concerning existence and uniqueness for SDE cannot be applied directly here for the coefficients of (3.3) are not Lipschitz. Furthermore, even if there exists a solution, one must show that the solution takes values in the set of states. Actually this property and the questions of existence and uniqueness are linked.

Concerning the property of being valued in the set of states, an important feature of the differential (3.3) is that it preserves the property to be a pure state (in quantum theory, a pure state is a one-dimensional projector). This idea is expressed in the following proposition.

PROPOSITION 3.1. *Let $(W_t)$ be a standard Brownian motion on $(\Omega, \mathcal{F}, \mathcal{F}_t, P)$ and let $|\psi_0\rangle$ be any norm-1 vector in $\mathbb{C}^2$. Let $\nu_t = \frac{1}{2}\langle\psi_t, (C + C^\star)\psi_t\rangle$ where $C$ is any operator on $\mathbb{C}^2$.*

*If the following stochastic equation*

$$(3.4) \quad d|\psi_t\rangle = (C - \nu_t I)|\psi_t\rangle\,dW_t + (-iH_0 - \tfrac{1}{2}(C^\star C - 2\nu_t C + \nu_t^2 I))|\psi_t\rangle\,dt$$

*admits a solution $(|\psi_t\rangle)$, then almost surely we have $\|\psi_t\| = 1$ for all $t$.*

*Furthermore the process $(|\psi_t\rangle\langle\psi_t|)$ takes values in the set of pure states and it is a solution of the diffusive Belavkin equation (3.3).*

PROOF. Let $|\psi_0\rangle$ be any vector in $\mathbb{C}^2$ and let $(|\psi_t\rangle)$ be a solution of (3.4). Let us prove that $\|\psi_t\|^2 = 1$. Using the Itô formulas and the fact that $H$ is self-adjoint, a straightforward computation shows that

$$\begin{aligned}
d\|\psi_t\|^2 &= d\langle\psi_t, \psi_t\rangle = \langle d\psi_t, \psi_t\rangle + \langle\psi_t, d\psi_t\rangle + \langle d\psi_t, d\psi_t\rangle \\
&= \langle(C - \nu_t I)\psi_t, \psi_t\rangle \, dW_t \\
&\quad + \langle(-iH_0 - \tfrac{1}{2}(C^\star C - 2\nu_t C + \nu_t^2 I))\psi_t, \psi_t\rangle \, dt \\
&\quad + \langle\psi_t, (C - \nu_t I)\psi_t\rangle \, dW_t \\
&\quad + \langle\psi_t, (-iH_0 - \tfrac{1}{2}(C^\star C - 2\nu_t C + \nu_t^2 I))\psi_t\rangle \, dt \\
&\quad + \langle(C - \nu_t I)\psi_t, (C - \nu_t I)\psi_t\rangle \, dt \\
&= (2\nu_t - 2\nu_t \langle\psi_t, \psi_t\rangle) \, dW_t.
\end{aligned}$$

If $\|\psi_0\|^2 = 1$, this implies that almost surely

$$\|\psi_t\|^2 = \|\psi_0\|^2 = 1$$

for all $t \geq 0$. Define the process $\rho_t = |\psi_t\rangle\langle\psi_t|$. It is valued in the set of pure states. As $\|\psi_t\| = 1$, we have for all $y \in \mathbb{C}^2$

$$\rho_t |y\rangle = \langle\psi_t, y\rangle |\psi_t\rangle.$$

Hence we can compute $d\rho_t |y\rangle$ by the Itô formula:

$$\begin{aligned}
d\rho_t |y\rangle &= \langle d\psi_t, y\rangle |\psi_t\rangle + \langle\psi_t, y\rangle d|\psi_t\rangle + \langle d\psi_t, y\rangle d|\psi_t\rangle \\
&= \langle(C - \nu_t)\psi_t, y\rangle |\psi_t\rangle \, dW_t \\
&\quad + \langle(-iH_0 - \tfrac{1}{2}(C^\star C - 2\nu_t C + \nu_t^2))\psi_t, y\rangle |\psi_t\rangle \, dt \\
&\quad + \langle\psi_t, y\rangle (C - \nu_t)|\psi_t\rangle \, dW_t \\
&\quad + \langle\psi_t, y\rangle (-iH_0 - \tfrac{1}{2}(C^\star C - 2\nu_t C + \nu_t^2))|\psi_t\rangle \, dt \\
&\quad + \langle(C - \nu_t)\psi_t, y\rangle (C - \nu_t)|\psi_t\rangle \, dt.
\end{aligned}$$

Let us show that this corresponds to (3.3). It is clear that $\nu_t = \tfrac{1}{2}\langle\psi_t, (C + C^\star)\psi_t\rangle$ corresponds to the term $\tfrac{1}{2}\operatorname{Tr}[|\psi_t\rangle\langle\psi_t|(C + C^\star)]$. As a consequence the term in front of the Brownian motion becomes

$$\begin{aligned}
&\langle(C - \nu_t)\psi_t, y\rangle |\psi_t\rangle + \langle\psi_t, y\rangle (C - \nu_t)|\psi_t\rangle \\
&= (C|\psi_t\rangle\langle\psi_t| + |\psi_t\rangle\langle\psi_t|C^\star - \operatorname{Tr}[|\psi_t\rangle\langle\psi_t|(C + C^\star)]|\psi_t\rangle\langle\psi_t|)|y\rangle.
\end{aligned}$$

A similar computation shows that the term in front of $dt$ is

$$L(|\psi_t\rangle\langle\psi_t|)|y\rangle.$$

Hence we recover the expression of Belavkin equation (3.3) and the proposition is proved. $\square$

As a consequence, we can express an existence and uniqueness theorem for (3.3). In what follows, we use the notion of "wave function." A wave function is a norm-1 vector which defines a pure state.

THEOREM 3.1.  *Let $(\Omega, \mathcal{F}, \mathcal{F}_t, P)$ be a probability space which supports a standard Brownian motion $(W_t)$ and let $\rho_0$ be any state on $\mathbb{C}^2$.*

*The stochastic differential equation*

$$\rho_t = \rho_0 + \int_0^t L(\rho_s)\, ds + \int_0^t [\rho_s C^\star + C\rho_s - \text{Tr}[(\rho_s(C + C^\star))\rho_s]]\, dW_s$$

*admits a unique solution $(\rho_t)$. The solution takes values in the set of states and is defined for all $t \geq 0$.*

*Furthermore, if the initial condition is a pure state, the solution takes values in the set of pure states. The corresponding stochastic differential equation for a wave function is then given by*

$$d|\psi_t\rangle = (C - \nu_t)|\psi_t\rangle\, dW_t + (-iH_0 - \tfrac{1}{2}(C^\star C - 2\nu_t C + \nu_t^2))|\psi_t\rangle\, dt$$

*where $\nu_t = \tfrac{1}{2}\langle \psi_t, (C + C^\star)\psi_t \rangle$.*

PROOF.  As the coefficients of (3.3) are not Lipschitz, we cannot apply directly the usual existence and uniqueness theorems for SDE. However, the coefficients are $C^\infty$, hence locally Lipschitz. We can use a truncature method. Equation (3.3) is of the following form:

(3.5) $$d\rho_t = L(\rho_t)\, dt + \Theta(\rho_t)\, dW_t$$

where $\Theta$ is $C^\infty$ and $\Theta(A) = AC^\star + CA - \text{Tr}[A(C + C^\star)]A$. Let $k > 1$ be an integer; we define the truncation function $\varphi_k$ from $\mathbb{R}$ to $\mathbb{R}$ defined by

$$\varphi_k(x) = \begin{cases} -k, & \text{if } x \leq -k, \\ x, & \text{if } -k \leq x \leq k, \\ k, & \text{if } -k \leq x \leq k. \end{cases}$$

For a matrix $A = (a_{ij})$, we define by extension $\tilde{\varphi}_k(A) = \varphi_k(\text{Re}(a_{ij})) + i\varphi_k(\text{Im}(a_{ij}))$. Thus the function $\Theta \circ \tilde{\varphi}_k$ is Lipschitz. Now we consider the truncated equation:

$$d\rho_{k,t} = L \circ \tilde{\varphi}_k(\rho_{k,t})\, dt + \Theta \circ \tilde{\varphi}_k(\rho_{k,t})\, dW_t.$$

The Cauchy–Lipschitz theorem concerning stochastic differential equations can be applied; there exists a unique solution $t \mapsto \rho_{k,t}$ defined for all $t$. Besides the solution is continuous in time.

Define the random stopping times

$$T_k = \inf\{t, \exists (ij) / |\text{Re}(a_{ij}(\rho_{k,t}))| = k \text{ or } |\text{Im}(a_{ij}(\rho_{k,t}))| = k\}.$$

As $\rho_0$ is a state, we have $\|\rho_0\| \leq 1$. Thanks to continuity, if $k$ is chosen large enough, we have $T_k > 0$ and for all $t \leq T_k$ we have $\tilde{\varphi}_k(\rho_{k,t}) = \rho_{k,t}$. Thus $t \mapsto \rho_{k,t}$ is the unique solution of (3.3) (without truncation) on $[0, T_k]$. The usual method for solving an equation with non-Lipschitz coefficients is to put $T = \lim_k T_k$ and to show that $T = \infty$.

In addition to the proof of existence of a solution, we must prove that the process is valued in the set of states. If $\nu$ is any state, we have $\|\nu\| \leq 1$, so $|\nu(ij)| \leq 1$. Hence if we prove that on $[0, T_2]$ the process $(\rho_{2,t})$ is valued on a set of states, this would prove that $T_2 = \infty$ a.s. and we would have proved that there exists a unique solution valued in the set of states. Let us prove this fact.

In the proof of the existence and uniqueness of a solution in the case of Cauchy–Lipschitz coefficients, the solution is obtained as the limit of the sequence

$$(3.6) \quad \begin{cases} \rho_{n+1}(t) = \rho_n(0) + \int_0^t L \circ \tilde{\varphi}_k(\rho_n(s))\, ds + \int_0^t \Theta \circ \tilde{\varphi}_k(\rho_n(s))\, dW_s, \\ \rho_0(t) = \rho. \end{cases}$$

With our definition of $\Theta$ and $L$, if $\rho_0$ is a state, it is clear that this sequence is self-adjoint with trace 1. These conditions are closed and at the limit the process is self-adjoint with trace 1. But the condition of positivity does not follow from such arguments. We shall prove it by other means.

Consider the random time

$$(3.7) \quad T^0 = \inf\{t \leq T_2 / \exists X \in \mathbf{C}^2 / \langle X, \rho_{2,t} X \rangle = 0\}.$$

We have $\langle X, \rho_0 X \rangle \geq 0$ for all $X$, so by continuity we have $\langle X, \rho_{2,t} X \rangle \geq 0$ on $[0, T^0]$ which implies that $\rho_{2,t}$ is a state for all $t \leq T^0$.

If $T^0 = T_2$ a.s., the result is proved. Otherwise, if we have $T^0 < T_2$, then by continuity there exists $X$ such that $\langle X, \rho_{2,T^0} X \rangle = 0$ and for all $Y$ $\langle Y, \rho_{2,T^0} Y \rangle \geq 0$. This implies that $\rho_{2,T^0}$ is a pure state because we work in dimension 2. Let us denote by $\psi_{T^0}$ a vector of norm 1 such that $\rho_{2,T^0} = |\psi_{T^0}\rangle\langle\psi_{T^0}|$. Consider the equation

$$d|\psi_t\rangle = (C - \nu_t)|\psi_t\rangle\, dW_t + (-iH_0 - \tfrac{1}{2}(C^\star C - 2\nu_t C + \nu_t^2))|\psi_t\rangle\, dt$$

with $\psi_{T^0}$ as initial condition. The problem of existence and uniqueness for this equation is solved by a truncation method also. The fact that, if we have a solution, it is of norm 1 shows that the solution obtained by truncation (defined for all $t$) is actually the solution of (3.4). Proposition 2.1 and the uniqueness of $\rho_{2,t}$ on $[T^0, T_2]$ show that the solution of (3.4) which satisfies

$$|\psi_t\rangle = |\psi_{T^0}\rangle + \int_{T^0}^t (C - \nu_s)|\psi_s\rangle\, dW_s + (-iH_0 - \tfrac{1}{2}(C^\star C - 2\nu_s C + \nu_s^2))|\psi_s\rangle\, ds$$

defines a process $(|\psi_t\rangle\langle\psi_t|)$ equal to $\rho_{2,t}$ on $[T^0, T_2]$. Hence the process obtained by truncation is valued on set of states and the theorem is proved. □

### 3.2. Change of measure.

At this stage, it must be said that the stochastic differential equation usually appearing in the literature is of the following form:

$$\rho_t = \rho_0 + \int_0^t L(\rho_s)\,ds + \int_0^t [\rho_s C^\star + C\rho_s - \mathrm{Tr}[\rho_s(C+C^\star)]]\,d\tilde{W}_s, \quad (3.8)$$

where

$$\tilde{W}_t = W_t - \int_0^t \mathrm{Tr}[\rho_s(C+C^\star)]\,ds. \quad (3.9)$$

Hence it seems to be rather different from (3.3). Actually the link between the two different equations is given by Girsanov's theorem (see [12]).

THEOREM 3.2. *Let $(W_t)$ be a standard Brownian motion on $(\Omega, \mathcal{F}, \mathcal{F}_t, P)$ and let $H$ be a càdlàg process. Let*

$$X_t = \int_0^t H_s\,ds + W_t \quad (3.10)$$

*and define a new probability by*

$$\frac{dQ}{dP} = \exp\left(-\int_0^T H_s\,dWs - \frac{1}{2}\int_0^T H_s^2\,ds\right)$$

*for some $T > 0$. Then under $Q$, the process $(X_t)$ is a standard Brownian motion for $0 \leq t \leq T$.*

The link between the two equations (3.3) and (3.8) is then obvious. Let $(\rho_t)$ be the solution of (3.3) given by Theorem 3.1 on $(\Omega, \mathcal{F}, \mathcal{F}_t, P)$. For some $T > 0$, define the probability measure $Q$ by

$$\frac{dQ}{dP} = \exp\left(\int_0^T \mathrm{Tr}[\rho_t(C+C^\star)]\,dW_s - \frac{1}{2}\int_0^T \mathrm{Tr}[\rho_t(C+C^\star)]^2\,ds\right). \quad (3.11)$$

The above theorem claims that $\tilde{W}_t$ is a standard Brownian motion under $Q$ for $0 \leq t \leq T$. Hence (3.8) is the same equation as (3.3) up to a change of measure. In the following section, we show that the solution of (3.3) can be obtained as limit of discrete quantum trajectories.

## 4. Approximation and convergence.

### 4.1. The discrete approximation.

In this section, we present a way to obtain the solution of the diffusive Belavkin equation (3.3) as a limit of concrete discrete quantum trajectories. Let us show that these discrete trajectories satisfy evolution equations which appear as approximations of stochastic differential equations.

In Section 1, we had the discrete quantum trajectories satisfying

$$(4.1) \quad \rho_{k+1} = \mathcal{L}_0(\rho_k) + \mathcal{L}_1(\rho_k) + \left[ -\sqrt{\frac{q_{k+1}}{p_{k+1}}}\mathcal{L}_0(\rho_k) + \sqrt{\frac{p_{k+1}}{q_{k+1}}}\mathcal{L}_1(\rho_k) \right] X_{k+1}.$$

Hence, we have

$$
\begin{aligned}
(4.2) \quad \rho_{k+1} - \rho_0 &= \sum_{i=0}^{k} [\rho_{i+1} - \rho_i] \\
&= \sum_{i=0}^{k} [\mathcal{L}_0(\rho_i) + \mathcal{L}_1(\rho_i) - \rho_i] \\
&\quad + \sum_{i=0}^{k} \left[ -\sqrt{\frac{q_{i+1}}{p_{i+1}}}\mathcal{L}_0(\rho_i) + \sqrt{\frac{p_{i+1}}{q_{i+1}}}\mathcal{L}_1(\rho_i) \right] X_{i+1}.
\end{aligned}
$$

Let us introduce a time discretization. Consider a partition of $[0,T]$ in subintervals of equal size $1/n$. The time of interaction is supposed now to be $h = 1/n$; the unitary operator depends then on the time interaction:

$$U(n) = \begin{pmatrix} L_{00}(n) & L_{01}(n) \\ L_{10}(n) & L_{11}(n) \end{pmatrix}.$$

In [1], Attal and Pautrat have shown that the asymptotic of the coefficients $L_{ij}(n)$ must be properly rescaled in order to obtain a nontrivial limit when $n$ goes to infinity. Indeed they have shown that $V_{[nt]} = U_{[nt]}(n)\ldots U_1(n)$, which represents the discrete dynamic of quantum repeated interactions, converges to a process $V_t$ solution of a quantum Langevin equation only if the coefficients $L_{ij}$ obey certain normalizations. When translated to our context, the results of [1] show that we should consider

$$(4.3) \qquad L_{00}(n) = I + \frac{1}{n}\left(-iH_0 - \frac{1}{2}CC^\star\right) + o\left(\frac{1}{n}\right),$$

$$(4.4) \qquad L_{10}(n) = \frac{1}{\sqrt{n}}C + o\left(\frac{1}{n}\right).$$

Remember that the unitary operator is given by

$$U(n) = \exp\left(i\frac{1}{n}H_{\text{tot}}\right).$$

The corresponding Hamiltonian $H_{\text{tot}}$ which gives the expression $U(n)$ is of the following form:

$$H_{\text{tot}} = H_0 \otimes I + I \otimes \begin{pmatrix} 1 & 0 \\ 0 & 0 \end{pmatrix} + \frac{1}{\sqrt{n}}\left[ C \otimes \begin{pmatrix} 0 & 0 \\ 1 & 0 \end{pmatrix} + C^\star \otimes \begin{pmatrix} 0 & 1 \\ 0 & 0 \end{pmatrix} \right] + o\left(\frac{1}{n}\right),$$

where $H_0$ is the Hamiltonian of the small system and $C$ is any operator on $\mathbb{C}^2$.

With the time discretization we then obtain

(4.5) $\rho_{k+1}(n) = \mathcal{L}_0(n)(\rho_k(n)) + \mathcal{L}_1(n)(\rho_k(n))$

$$+ \left[-\sqrt{\frac{q_{k+1}(n)}{p_{k+1}(n)}}\mathcal{L}_0(n)(\rho_k(n)) + \sqrt{\frac{p_{k+1}(n)}{q_{k+1}(n)}}\mathcal{L}_1(n)(\rho_k(n))\right]$$

(4.6)
$$\times X_{k+1}(n).$$

Remember that the sequence of random variables $(X_k(n))$ is defined through the two probabilities:

(4.7) $\begin{cases} p_{k+1} = \text{Tr}[\mathcal{L}_0(\rho_k)], \\ q_{k+1} = \text{Tr}[\mathcal{L}_1(\rho_k)]. \end{cases}$

By definition of $(X_k)$ we have

(4.8) $X_k(n)(i) = \begin{cases} -\sqrt{\dfrac{q_{k+1}(n)}{p_{k+1}(n)}}, & \text{with probability } p_{k+1}(n) \text{ if } i=0, \\ \sqrt{\dfrac{p_{k+1}(n)}{q_{k+1}(n)}}, & \text{with probability } q_{k+1}(n) \text{ if } i=1. \end{cases}$

Each probability depends on the expression of $\mathcal{L}_i$, which depends on the measured observable $A = \lambda_0 P_0 + \lambda_1 P_1$. At the limit, the diffusive behavior of $\rho_{[nt]}(n)$ is depending on the comportment of $(X_k)$, that is, on the observable:

1. If the observable is of the form $A = \lambda_0 \begin{pmatrix} 1 & 0 \\ 0 & 0 \end{pmatrix} + \lambda_1 \begin{pmatrix} 0 & 0 \\ 0 & 1 \end{pmatrix}$, we obtain the following asymptotic for the probabilities:

$$p_{k+1}(n) = 1 - \frac{1}{n}\text{Tr}[\mathcal{J}(\rho_k(n))] + o\left(\frac{1}{n}\right),$$

$$q_{k+1}(n) = \frac{1}{n}\text{Tr}[\mathcal{J}(\rho_k(n))] + o\left(\frac{1}{n}\right).$$

The discrete equation then becomes

$$\rho_{k+1}(n) - \rho_k(n)$$
$$= \frac{1}{n}L(\rho_k(n)) + o\left(\frac{1}{n}\right) + \left[\frac{\mathcal{J}(\rho_k(n))}{\text{Tr}[\mathcal{J}(\rho_k(n))]} - \rho_k(n) + o(1)\right]$$
$$\times \sqrt{q_{k+1}(n)p_{k+1}(n)}X_{k+1}(n).$$

2. If the observable $A$ is nondiagonal in the basis $(\Omega, X)$, for the eigenprojectors, put $P_0 = \begin{pmatrix} p_{00} & p_{01} \\ p_{10} & p_{11} \end{pmatrix}$ and $P_1 = \begin{pmatrix} q_{00} & q_{01} \\ q_{10} & q_{11} \end{pmatrix}$; we have

$$p_{k+1} = p_{00} + \frac{1}{\sqrt{n}}\text{Tr}[\rho_k(p_{01}C + p_{10}C^\star)] + \frac{1}{n}\text{Tr}[\rho_k p_{00}(C+C^\star)] + o\left(\frac{1}{n}\right),$$

$$q_{k+1} = q_{00} + \frac{1}{\sqrt{n}} \operatorname{Tr}[\rho_k(q_{01}C + q_{10}C^\star)] + \frac{1}{n} \operatorname{Tr}[\rho_k q_{00}(C + C^\star)] + o\left(\frac{1}{n}\right).$$

The discrete equation here becomes

(4.9) $$\rho_{k+1} - \rho_k = \frac{1}{n}L(\rho_k) + o\left(\frac{1}{n}\right)$$

(4.10)
$$+ [e^{i\theta}C\rho_k + e^{-i\theta}\rho_k C^\star$$
$$- \operatorname{Tr}[\rho_k(e^{i\theta}C + e^{-i\theta}C^\star)]\rho^k + o(1)]\frac{1}{\sqrt{n}}X_{k+1}.$$

In this expression, the parameter $\theta$ represents a kind of phase. It is real and depends on the coefficients of the eigenprojectors. If we put $C_\theta = e^{i\theta}C$, the discrete (4.9) becomes

$$\rho_{k+1} - \rho_k = \frac{1}{n}L(\rho_k) + o\left(\frac{1}{n}\right)$$
$$+ [C_\theta\rho_k + \rho_k C_\theta^\star - \operatorname{Tr}[\rho_k(C_\theta + C_\theta^\star)]\rho_k + o(1)]\frac{1}{\sqrt{n}}X_{k+1}.$$

For each $\theta$, we have similar expressions for discrete equations with different operators $C_\theta$. Let us stress that this parameter does not modify the expression of $L$. In the following, we deal with $\theta = 0$.

In [11], it is shown that the case where $A$ is diagonal gives rise to the jump-Belavkin equation at a continuous limit. In the following section, we show that the diffusive case is obtained as the limit of the discrete process which comes from the measurement of a nondiagonal observable.

4.2. *Convergence theorems.* Before presenting the main theorem concerning the convergence of discrete quantum trajectories, we show a first result concerning the average of the processes. In order to avoid confusion between the discrete-time process $(\rho_k)$ and the continuous-time process $(\rho_t)$ we write the discrete process $(\rho^k)$ with the index on the top.

THEOREM 4.1. *Let $(\Omega, X)$ be any orthonormal basis of $\mathbb{C}^2$. For all nondiagonal observable $A$, the deterministic function $t \mapsto \to \mathbf{E}[\rho^{[nt]}(n)]$ converges in $L^\infty([0,T])$, when $n$ goes to infinity, to the function $t \mapsto \to \mathbf{E}[\rho_t]$. That is,*

$$\sup_{0 < s < T} \|\mathbf{E}[\rho^{[ns]}(n)] - \mathbf{E}[\rho_s]\| \underset{n \to \infty}{\longrightarrow} 0.$$

*Furthermore the function $t \mapsto \to \mathbf{E}[\rho_t]$ is the solution of the Master equation*

$$d\nu_t = L(\nu_t)\,dt.$$

PROOF. First of all, we show the second part of the theorem. We can consider the function $t \mapsto \to \mathbf{E}[\rho_t]$ because we have existence and uniqueness of the solution of (3.3). The process $(\rho_t)$ is integrable (because $\rho_t$ is a state for all $t$). It is obvious that this function takes also values in the set of states. As $\rho_0$ is a deterministic state we must show

$$\mathbf{E}[\rho_t] = \rho_0 + \int_0^t L(\mathbf{E}[\rho_s])\,ds. \tag{4.11}$$

We know that the process $(\rho_t)$ satisfies

$$\rho_t = \rho_0 + \int_0^t L(\rho_s)\,ds + \int_0^t [\rho_s C^\star + C\rho_s - \mathrm{Tr}(\rho_s(C+C^\star))\rho_s]\,dW_s.$$

As $(W_t)$ is a martingale and the process $(\rho_t)$ is predictable (for it is continuous), the properties of stochastic integral give

$$\mathbf{E}\left[\int_0^t [\rho_s C^\star + C\rho_s - \mathrm{Tr}(\rho_s(C+C^\star))\rho_s]\,dW_s\right] = 0.$$

Hence, we have, by linearity of $L$,

$$\mathbf{E}[\rho_t] = \rho_0 + \int_0^t \mathbf{E}[L(\rho_s)]\,ds$$
$$= \rho_0 + \int_0^t L(\mathbf{E}[\rho_s])\,ds.$$

We then have the integral form of the solution of the Master equation and the second part is proved.

Let us show now the first part of the theorem. We shall now compare $\mathbf{E}[\rho^{[nt]}(n)]$ with $\mathbf{E}[\rho_t]$. As in the continuous case, the martingale argument is replaced by the fact that the process $(X_k)$ is centered. Remember that we have

$$\mathbf{E}[X_{k+1}] = \mathbf{E}[\mathbf{E}[X_{k+1}/\mathcal{F}_k]] = 0.$$

As a consequence, considering $k = [nt]$ and taking expectation in the discrete equation, we have

$$\mathbf{E}[\rho^{[nt]}(n)] - \rho^0 = \sum_{i=0}^{[nt]-1} \frac{1}{n} L(\mathbf{E}[\rho^k(n)]) + o\left(\frac{1}{n}\right).$$

This is a kind of Euler scheme and we can conclude by a discrete Gronwall lemma argument that we have

$$\sup_{0 < s < t} \|\mathbf{E}[\rho^{[ns]}(n)] - \mathbf{E}[\rho_s]\| \xrightarrow[n \to \infty]{} 0. \qquad \square$$

The average of the discrete process is then an approximation of the average of $\rho_t$. In [1], this result was shown in the case of repeated interactions without measurement. It is a consequence of the asymptotic of the unitary-operator coefficients, so it justifies our choice of the coefficients.

Concerning the convergence of the processes, the discrete process which is the candidate to converge to the diffusive quantum trajectory satisfies for $k = [nt]$

$$\rho^{[nt]} - \rho^0 = \sum_{i=0}^{[nt]-1} \frac{1}{n} L(\rho^k(n)) + o\left(\frac{1}{n}\right) + \sum_{i=0}^{[nt]-1} [\Theta(\rho^k) + o(1)] \frac{1}{\sqrt{n}} X_{i+1}.$$

Thanks to this equation we can define the processes:

(4.12)
$$W_n(t) = \frac{1}{\sqrt{n}} \sum_{k=1}^{[nt]} X_k(n),$$

$$V_n(t) = \frac{[nt]}{n},$$

$$\rho_n(t) = \rho^{[nt]}(n),$$

$$\varepsilon_n(t) = \sum_{i=0}^{[nt]-1} o\left(\frac{1}{n}\right) + \sum_{i=0}^{[nt]-1} o(1) \frac{1}{\sqrt{n}} X_{i+1}.$$

By observing that these four processes are piecewise constant, we can write the process $(\rho_n(t))_{t \geq 0}$ as a solution of a stochastic differential equation in the following way:

(4.13) $$\rho_n(t) = \rho_0 + \varepsilon_n(t) + \int_0^t L(\rho_n(s-))\,dV_n(s) + \int_0^t \Theta(\rho_n(s-))\,dW_n(s).$$

We now use a theorem of Kurtz and Protter (cf. [8]) to prove the convergence. Let us first fix some notation.

Recall that the square-bracket $[X, X]$ is defined for a semimartingale by the formula

$$[X, X]_t = X_t^2 - 2 \int_0^t X_{s-}\,dX_s.$$

We shall denote by $T_t(V)$ the total variation of a finite-variation process $V$ on the interval $[0, t]$. The theorem of Kurtz and Protter [8] (see [7]) that we use is the following.

THEOREM 4.2. *Let $W_n$ be a martingale and let $V_n$ be a finite variation process. Consider the process $X_n$ defined by*

$$X_n(t) = \rho_0 + \varepsilon_n(t) + \int_0^t L(X_n(s-))\,dV_n(s) + \int_0^t \Theta(X_n(s-))\,dW_n(s).$$

*Assume that for each $t \geq 0$,*

$$\sup_n \mathbf{E}[[W_n, W_n]_t] < \infty,$$

$$\sup_n \mathbf{E}[T_t(V_n)] < \infty,$$

*and that $(W_n, V_n, \varepsilon_n)$ converges in distribution to $(W, V, 0)$ where $W$ is a standard Brownian motion and $V(t) = t$ for all $t$. Suppose that $X$ satisfies*

$$X_t = X_0 + \int_0^t L(X_s)\,ds + \int_0^t \Theta(X_s)\,dW_s$$

*and that the solution of this stochastic differential equation is unique. Then $X_n$ converges in distribution to $X$.*

We wish to apply this theorem to the process $(\rho_n(t))$ [(4.13)]. The first step is the convergence of $W_n$ in (4.12) to a Brownian motion. We need the following theorem (cf. [3, 9]) which is a generalization of the Donsker theorem.

THEOREM 4.3. *Let $(M_n)$ be a sequence of martingales. Suppose that*

$$\lim_{n \to \infty} \mathbf{E}\left[\sup_{s \leq t} |M_n(s) - M_n(s-)|\right] = 0$$

*and*

$$[M_n, M_n]_t \underset{n \to \infty}{\longrightarrow} t.$$

*Then $M_n$ converges in distribution to a standard Brownian motion. The conclusion is the same if we just have*

$$\lim_{n \to \infty} \mathbf{E}[|[M_n, M_n]_t - t|] = 0.$$

Back to our process $W_n$, consider the filtration

$$\mathcal{F}_t^n = \sigma(X_i, i \leq [nt]).$$

PROPOSITION 4.1. *Let $(W_n, V_n, \varepsilon_n)$ be the processes defined in (4.12). We have that $(W_n(t))$ is a $\mathcal{F}_t^n$-martingale. The process $(W_n)$ converges to a standard Brownian motion $W$ when $n$ goes to infinity. Furthermore we have*

$$\sup_n \mathbf{E}[[W_n, W_n]_t] < \infty.$$

*Finally, we have the convergence in distribution for the processes $(W_n, V_n, \varepsilon_n)$, when $n$ goes to infinity, to $(W, V, 0)$.*

PROOF. Thanks to the definition of the random variable $X_k$, we have $\mathbf{E}[X_{i+1}/\mathcal{F}_i^n]=0$ which implies $\mathbf{E}[\frac{1}{n}\sum_{i=[ns]+1}^{[nt]} X_i/\mathcal{F}_s^n]=0$ for $t>s$. Thus if $t>s$, we have the martingale property

$$\mathbf{E}[W_n(t)/\mathcal{F}_s^n] = W_n(s) + \mathbf{E}\left[\frac{1}{\sqrt{n}} \sum_{i=[ns]+1}^{[nt]} X_i/\mathcal{F}_s^n\right] = W_n(s).$$

By definition of $[Y,Y]$ for a stochastic process we have

$$[W_n,W_n]_t = W_n(t)^2 - 2\int_0^t W_n(s-)\,dW_n(s) = \frac{1}{n}\sum_{i=1}^{[nt]} X_i^2.$$

Thus we have

$$\mathbf{E}[[W_n,W_n]_t] = \frac{1}{n}\sum_{i=1}^{[nt]} \mathbf{E}[X_i^2] = \frac{1}{n}\sum_{i=1}^{[nt]} \mathbf{E}[\mathbf{E}[X_i^2/\sigma\{X_l, l<i\}]]$$

$$= \frac{1}{n}\sum_{i=1}^{[nt]} 1 = \frac{[nt]}{n}.$$

Hence we have

$$\sup_n \mathbf{E}[[W_n,W_n]_t] \leq t < \infty.$$

Let us prove the convergence of $(W_n)$ to $(W)$. According to Theorem 4.3 we must prove that

$$\lim_{n\to\infty} \mathbf{E}[|[M_n,M_n]_t - t|] = 0.$$

Actually we prove an $L_2$ convergence:

$$\lim_{n\to\infty} \mathbf{E}[|[M_n,M_n]_t - t|^2] = 0,$$

which implies the $L_1$ convergence. In order to show this convergence, we use the following property:

$$\mathbf{E}[X_i^2] = \mathbf{E}[\mathbf{E}[X_i^2/\sigma\{X_l,l<i\}]] = 1,$$

and if $i<j$,

$$\mathbf{E}[(X_i^2-1)(X_j^2-1)] = \mathbf{E}[\mathbf{E}[(X_i^2-1)(X_j^2-1)/\sigma\{X_l,l<j\}]]$$

$$= \mathbf{E}[(X_i^2-1)]\mathbf{E}[(X_j^2-1)]$$

$$= 0.$$

This gives

$$\mathbf{E}\bigg[\bigg([W_n, W_n]_t - \frac{[nt]}{n}\bigg)^2\bigg] = \frac{1}{n^2}\sum_{i=1}^{[nt]} \mathbf{E}[(X_i^2 - 1)^2] + \frac{1}{n^2}\sum_{i<j}\mathbf{E}[(X_i^2 - 1)(X_j^2 - 1)]$$

$$= \frac{1}{n^2}\sum_{i=1}^{[nt]} \mathbf{E}[(X_i^2 - 1)^2].$$

Thanks to the fact that $p_{00}$ and $q_{00}$ are not equal to zero (because the observable $A$ is not diagonal!) the terms $\mathbf{E}[(X_i^2 - 1)^2]$ are bounded uniformly in $i$ so we have

$$\lim_{n\to\infty} \mathbf{E}\bigg[\bigg([W_n, W_n]_t - \frac{[nt]}{n}\bigg)^2\bigg] = 0.$$

As $\frac{[nt]}{n} \longrightarrow t$ in $L_2$ we have the desired convergence. Finally, the convergence in distribution of $(W_n)$ and $(V_n)$ implies the convergence of $(\varepsilon_n)$ to 0.  □

We can now express the final theorem.

THEOREM 4.4.  *Let $(\Omega, X)$ be any orthonormal basis of $\mathbb{C}^2$ and let $A$ be any nondiagonal observable (in this basis). Let $\rho_0$ be any initial state on $\mathbb{C}^2$. Let $(\rho^{[nt]}(n))$ be the discrete quantum trajectory obtained from the quantum repeated measurement principle with respect to $A$. The process $(\rho^{[nt]}(n))$ then satisfies*

$$\rho^{[nt]}(n) = \rho_0 + \sum_{i=0}^{[nt]-1} \frac{1}{n} L(\rho^k(n)) + o\bigg(\frac{1}{n}\bigg) + \sum_{i=0}^{[nt]-1} [\Theta(\rho^k) + o(1)]\frac{1}{\sqrt{n}} X_{i+1}.$$

*Let $(\rho_t)$ be the solution of the diffusive Belavkin equation (3.3) which satisfies*

$$\rho_t = \rho_0 + \int_0^t L(\rho_s)\, ds + \int_0^t \Theta(\rho_s)\, dW_s.$$

*Then we have the convergence in distribution*

$$(\rho^{[nt]}(n)) \xrightarrow[n\to\infty]{\mathcal{D}} (\rho_t).$$

PROOF.  It is a simple compilation of Theorems 4.2, Proposition 4.1 and existence and uniqueness of Theorem 3.1.  □

## REFERENCES

[1] ATTAL, S. and PAUTRAT, Y. (2006). From repeated to continuous quantum interactions. *Ann. Henri Poincaré* **7** 59–104. MR2205464

[2] BOUTEN, L., GUȚĂ, M. and MAASSEN, H. (2004). Stochastic Schrödinger equations. *J. Phys. A* **37** 3189–3209. [MR2042615](MR2042615)
[3] DACUNHA-CASTELLE, D. and DUFLO, M. (1986). *Probability and Statistics*. **II**. Springer, New York.
[4] DAVIES, E. B. (1976). *Quantum Theory of Open Systems*. Academic Press, London. [MR0489429](MR0489429)
[5] GOUGH, J. and SOBOLEV, A. (2004). Stochastic Schrödinger equations as limit of discrete filtering. *Open Syst. Inf. Dyn.* **11** 235–255. [MR2102379](MR2102379)
[6] KÜMMERER, B. and MAASSEN, H. (2003). An ergodic theorem for quantum counting processes. *J. Phys. A* **36** 2155–2161. [MR1963955](MR1963955)
[7] KURTZ, T. G. and PROTTER, P. (1991). Weak limit theorems for stochastic integrals and stochastic differential equations. *Ann. Probab.* **19** 1035–1070. [MR1112406](MR1112406)
[8] KURTZ, T. G. and PROTTER, P. (1991). Wong–Zakai corrections, random evolutions, and simulation schemes for SDEs. In *Stochastic Analysis* 331–346. Academic Press, Boston, MA. [MR1119837](MR1119837)
[9] KURTZ, T. G. and PROTTER, P. E. (1996). Weak convergence of stochastic integrals and differential equations. In *Probabilistic Models for Nonlinear Partial Differential Equations (Montecatini Terme, 1995). Lecture Notes in Math.* **1627** 1–41. Springer, Berlin. [MR1431298](MR1431298)
[10] PARTHASARATHY, K. R. (1992). *An Introduction to Quantum Stochastic Calculus. Monographs in Mathematics* **85**. Birkhäuser, Basel. [MR1164866](MR1164866)
[11] PELLEGRINI, C. (2007). Existence, uniqueness and approximation of stochastic Schrödinger equation: the Poissonian case. To appear.
[12] PROTTER, P. E. (2004). *Stochastic Integration and Differential Equations*, 2nd ed. *Applications of Mathematics (New York)* **21**. Springer, Berlin. [MR2020294](MR2020294)

INSTITUT C. JORDAN
UNIVERSITÉ C. BERNARD, LYON 1
21 AV CLAUDE BERNARD
69622 VILLEURBANNE CEDEX
FRANCE
E-MAIL: [pelleg@math.univ-lyon1.fr](pelleg@math.univ-lyon1.fr)